\documentclass[reqno]{amsart}   
\usepackage{pb-diagram,enumerate,amsmath}   
 \usepackage{amsfonts,amssymb,amsmath,a4wide} 
\usepackage{pstricks,pst-node,pst-coil} 
   
\theoremstyle{plain}    
\newtheorem{thm}{Theorem}[section]   
\newtheorem{cor}[thm]{Corollary}   
\newtheorem{lem}[thm]{Lemma}   
\newtheorem{prop}[thm]{Proposition}   
\theoremstyle{remark}    
\newtheorem{rem}[thm]{Remark}   
\theoremstyle{definition}    
   
\def\al{{\alpha}}       
\def\be{{\beta}}       
\def\de{{\delta}}       
       
\def\om{{\omega}}       
       
\def\la{{\lambda}}

\def\si{{\sigma}}       
       
\def\ga{{\gamma}}       
\def\ep{{\varepsilon}}

\def\phi{{\varphi}}

\let\pa\partial   
\let\na\nabla   
   
\DeclareMathAlphabet{\doba}{U}{msb}{m}{n}

\gdef\mR{\doba{R}}       
\gdef\mS{\doba{S}}

\def\cJ{{\mathcal J}}

\def\lamin{\lambda_{\rm min}^+}   
 
\def\Vol{{\mathop{\rm Vol}}}   
\def\Scal{{\mathop{\rm Scal}}}   
\def\Ric{{\mathop{\rm Ric}}}   
\def\ijkdiff{{i,j,k\atop i\neq j \neq k \neq i}}   
\def\ijkabc{ {ijk \al \be \ga \atop i\neq j \neq k \neq i}} 
\let\ti\tilde 
\def\eref#1{{\rm (\ref{#1})}} 
\def\Spin{\mathrm{Spin}}
\def\SO{\mathrm{SO}}
\def\geucl{g_{\rm eucl}}

\let\ol\overline

\long\def\komment#1{}
 
\begin{document}   
\title{A spinorial analogue of Aubin's inequality}   
\author{B. Ammann, J.-F. Grosjean, E. Humbert and B. Morel}   
\date{March 2007}
   
\begin{abstract}   
Let $(M,g,\si)$ be a compact Riemannian spin manifold of dimension $\geq 2$. 
For any metric
$\tilde g$ conformal to $g$, we denote by  $\tilde\lambda$ the 
first positive eigenvalue of the Dirac operator on 
$(M,\tilde g,\si)$. 
We show that 
$$\inf_{\tilde{g} \in [g]} \tilde\lambda\; \Vol(M,\tilde g)^{1/n} \leq (n/2)\; \Vol(S^n)^{1/n}.$$
This inequality is a spinorial analogue of Aubin's inequality, an important 
inequality in the solution of the Yamabe problem.
The inequality is already known in the case $n \geq 3$ and in the 
case $n = 2$, 
$\ker D=\{0\}$. Our proof also works in the remaining case $n=2$, 
$\ker D\neq \{0\}$. With the same method we also prove that any conformal class
on a Riemann surface contains a metric with $2\tilde\lambda^2\leq \tilde\mu$, 
where
$\tilde\mu$ denotes the first positive eigenvalue of the Laplace operator.
\end{abstract}

\maketitle   

\bigskip
\begin{center}Contents\end{center}

\contentsline {section}{\tocsection {}{1}{Introduction}}{2}
\contentsline {section}{\tocsection {}{2}{A variational formulation fo the spin conformal invariant}}{3}
\contentsline {section}{\tocsection {}{3}{The Bourguignon-Gauduchon-trivialization}}{3}
\contentsline {section}{\tocsection {}{4}{Development of the metric at the point $p$}}{6}
\contentsline {section}{\tocsection {}{5}{The test spinor }}{8}
\contentsline {section}{\tocsection {}{6}{The proof of Theorem 1.1\hbox {} for $n\geq 3$}}{9}
\contentsline {section}{\tocsection {}{7}{The case $n=2$}}{17}
\contentsline {section}{\tocsection {}{}{References}}{21}
  
   
\section{Introduction}   
 
Let $(M,g)$ be a compact Riemannian manifold of dimension $n\geq 2$. We assume that $M$ is spin, and we fix a spin structure $\si$ on $M$.
For any metric $\tilde{g}$ in the conformal class $[g]$ of $g$, we   
write $\lambda_1^+(\tilde{g})$ for the smallest positive eigenvalue of the 
Dirac operator with respect to $(M,\tilde g,\si)$. We define     
$$\lamin(M,g,\si) = \inf_{\tilde{g} \in [g] } \lambda_1^+(\tilde{g})   
\Vol(M,\tilde{g})^{1/n}.$$
If $(M,g)$ is the round sphere $\mS^n$ 
equipped with the unique spin structure on 
$\mS^n$, we simply write $\lamin(\mS^n)$.    
It was proven in \cite{lott:86} ($\ker D=\{0\}$) and \cite{ammann:03} 
($\ker D\neq \{0\}$) that    
$$\lamin(M,g,\si) >0.$$   
Several articles have been devoted to the study of this spin-conformal   
invariant. 
A non-exhaustive list is
\cite{hijazi:86,lott:86,baer:92b,ammann.habil}. In this article we will prove
the following.
\begin{thm}\label{main}
Let $(M,g,\si)$ be a compact spin manifolds of dimension $n \geq 2$.      
Then,  
\begin{equation}\label{ineq.large}   
  \lamin(M,g,\si) \leq \lamin(\mS^n) = \frac{n}{2}\, \om_n^{\frac{1}{n}} 
\end{equation} 
where $\om_n$ stands for the volume of the standard sphere   
$\mS^n$.
\end{thm}

The proof relies on constructing a suitable test spinor. The main idea of 
this construction is to start with a Killing spinor on the round sphere. 
Under stereographic projection this spinor then yields a solution
to the equation $D\psi= c |\psi|^{2/(n-1)}\psi$ on flat $\mR^n$. This solution 
will be rescaled, cut off and finally transplanted to a neighborhood of 
a given point $p$ of the manifold $M$. For this transplantation we carry out
several calculations in a well-adapted trivialization of the spinor bundle.

The first steps in our proof are common in all dimensions.  
However, in some final estimates 
one has to distinguish between the cases $n\geq 3$ and $n=2$.

In dimension $n\geq 3$ two other proofs for the theorem have already 
been published: a geometric construction \cite[Theorem~3.1]{ammann:03} and a 
proof using an invariant for non-compact spin manifolds \cite{grosse:06}.
In these dimensions, it is mostly the method of proof that is interesting 
and helpful: the trivialization presented here has less terms in the
Taylor expansion than the trivialization 
by using parallel transport along radial geodesics.
Some formulae of our article also enter
in \cite{grosse:06}. The calculations of our article also provide helpful 
formulae used in \cite{ammann.humbert.morel:p03b}, \cite{ammann.humbert:03} and
\cite{raulot:06}.

The main interest of the theorem however lies in the case $n=2$.
The easier subcase $n=2$, $\ker D=\{0\}$ could be dealt with 
by a modification 
of the geometric proof \cite[Theorem~3.2]{ammann:03}, but the subcase
$n=2$ and $\ker D\neq \{0\}$ remained open for longtimes. Gro\ss{}e's method 
fails as well for $n=2$ as the contribution of a cut-off function in 
\cite[Lemma 2.1(ii)]{grosse:06} is too large. We assume that her method can be
adapted by using a logarithmic cut-off function, but the details have not 
been worked out yet.

Our method of proof in dimension $2$ actually admits applications to other
problems as well. For example, one obtains the following proposition
that provides a negative answer to a question raised in \cite{aaf:99}.

\begin{prop}[See Corollary \ref{compar}]\label{prop.compar}
Let $(M,g)$ be a Riemann surface with fixed spin structure $\sigma$. 
For any metric
$\bar{g}$ in the conformal class $[g]$, let $\mu_1(\bar{g})$ be the first
positive eigenvalue of the Laplacian, and let $\la^+_1(\bar{g})$ be the 
first positive eigenvalue of the Dirac operator on $(M,\bar{g},\sigma)$.
Then
$$\inf_{\bar{g}\in[g]} \frac{\la^+_1(\bar{g})^2 }{\mu_1(\bar{g})} \leqslant  \frac{1}{2}.$$
\end{prop}

\noindent The article is organized as follows: 
  in Section~\ref{funct1},   
we recall that $\lamin(M,g,\si)$   
has a variational characterization. Then, in Section~\ref{ident}  
we introduce a well-adapted local trivialization of the spinor bundle,
called the Bourguignon-Gauduchon-trivialization. 
In Section~\ref{dev} we calculate the first terms of the Taylor    
development of the Dirac operator in this trivialization. 
In the following, i.e.\ in Section~\ref{testspinor},  
we construct a good test spinor using a   
Killing spinor on $\mS^n$, and then in  Section~\ref{estimate}, we 
set this spinor in the functional to get Theorem~\ref{main} 
in dimension $n \geq 3$. In the last section, i.e.\ in Section~\ref{dim2},
we describe the modifications for the case $n=2$ and prove the proposition.

\noindent  {\it Acknowledgments}. 
The authors want to thank Oussama Hijazi (Nancy) for his support 
and encouragement for working at this article. 
B.~Ammann wants to thank C. B\"ar for some discussions about related subjects.
We thank the referee for drawing our attention to
the article \cite{takahashi:02}. B.~Ammann thanks the Max-Planck institute
for gravitational physics, at Potsdam-Golm for its hospitality.


\section{A variational formulation fo the spin conformal invariant} \label{funct1}    
For a section $\psi\in\Gamma(\Sigma M)$ we define    
   
$$J(\psi)=\frac{\Big(\int_M|D\psi|^{\frac{2n}{n+1}}v_g\Big)^\frac{n+1}{n}}{\int_M \langle D\psi,\psi\rangle v_g}.$$  
At some places we will wirte $J_g$ instead of $J$ inorder to indicate, that the functional is defined with respect to $g$.
Based on some idea from \cite{lott:86},
Ammann proved in \cite{ammann.habil} that    
\begin{eqnarray} \label{funct}   
\lamin(M,g,\si)=\inf_\psi J(\psi)   
\end{eqnarray}   
where the infimum is taken over the set of smooth spinor fields for which    
$$\left(\int_M \langle D\psi,\psi\rangle v_g \right)>0.$$   
Hence, to prove Theorem \ref{main},   
we are reduced to find a smooth spinor field $\psi$ satisfying the   
condition below and such that $J(\psi) \leq  \lamin(\mS^n)+ \ep$ where
$\ep>0$ is  arbitrary small.


\section{The Bourguignon-Gauduchon-trivialization}  \label{ident}   
 
As already explained before, 
the proof of our main theorem is based on a the construction 
of a suitable test spinor. We
first construct a ``good'' spinor field of $\mR^n$ 
and then transpose it on the manifold. 
In order to carry this out, 
we need to locally identify spinor fields on $(\mR^n,\geucl)$ 
and spinor fields on $(M,g)$. Such an identification will be provided 
by a well-adapted local trivialization
of the spinor bundle of $\Sigma(M,g)$. 

If a spin manifold $N$ 
carries two metrics $g_1$ and $g_2$, then it is a priori 
unclear how to identify spinors on $(N,g_1)$ and spinors on $(N,g_2)$.
Bourguignon and Gauduchon \cite{bourguignon.gauduchon:92} constructed a 
convenient map from the spinor bundle of $(N,g_1)$ to the spinor bundle
of $(N,g_2)$ that allows us to identify spinors, and it is this identification
that will provide the necessary identification to us.
The trivialization
will be called \emph{Bourguignon-Gauduchon-trivialization}.

This trivialization is more efficient than the commonly used 
``trivialization by
parallel transport along radial geodesics'': with respect to the
Bourguignon-Gauduchon-trivialization less terms appear in the Taylor expansion 
in Section~\ref{dev}. 
   
Let $(M,g)$ be a Riemannian manifold with a spin structure
$\sigma:\Spin(M,g)\to \SO(M,g)$.  
Let $(x_1,\ldots x_n)$ be the Riemannian normal coordinates given by the    
exponential map at $p\in M$:   
\begin{eqnarray*}   
\exp_p : U\subset T_pM\cong\mR^n  &\longrightarrow& V\subset M\\   
(x_1,\ldots,x_n)&\longmapsto&m   
\end{eqnarray*}   

Let    
\begin{eqnarray*}   
 G: V&\longrightarrow&S^2_+(n,\mR)\\   
m&\longmapsto& G_m:=(g_{ij}(m))_{ij}    
\end{eqnarray*}   
denote the smooth map which associates to any point $m\in V$, the matrix of   
the coefficients of the metric~$g$ at this point, expressed in the basis   
$(\partial_i:=\frac{\partial}{\partial x^i})_{1\leq i\leq n}\;$.  
Since $G_m$ is    
symmetric and positive definite, there is a unique symmetric and positive
definite matrix  
$B_m$ such that $$B_m^2=G_m^{-1}\;.$$  
Since $${}^t(B_mX)G_m(B_m Y)=\geucl(X,Y)\;,\quad\forall X,Y\in\mR^n\;,$$  
where $\geucl$ stands for the Euclidean    
scalar product, we get the following isometry defined by   
\begin{eqnarray*}   
B_m : (T_{\exp_p^{-1}(m)}U\cong\mR^n,\geucl)  &\longrightarrow& (T_mV,g_m)\\   
(a^1,\ldots,a^n)&\longmapsto&\sum_{i,j} b_i^j(m) a^i\partial_j(m)   
\end{eqnarray*}   
for each point $m\in V$, where $b_i^j(m)$ are the coefficients of the matrix $B_m$ (from now on, we use Einstein's summation convention).  
As the matrix $B_m$ depends smoothly on    
$m$, we can identify the following $\mathrm{SO}_n$-principal bundles:   
   
\dgARROWLENGTH=2em   
$$\begin{diagram}   
\node{\mathrm{SO}(U,\geucl)}\arrow[2]{e,t}{\eta}\arrow{s}\node[2]{\mathrm{SO}(V,g)}\arrow{s}\\   
\node{U\subset T_pM}\arrow[2]{e,t}{\exp_p}\node[2]{V\subset M}   
\end{diagram}$$   
where $\eta$ is given by the action of $B$ on each component vector
of a frame in $\mathrm{SO}(U,\geucl)$.    
The map $\eta$ commutes with the right action of $\mathrm{SO}_n$, 
therefore the map $\eta$ can be lifted to the spin structures
   
\dgARROWLENGTH=2em   
$$\begin{diagram}   
\node{\mathrm{Spin}(n)\times U=\mathrm{Spin}(U,\geucl)}\arrow[2]{e,t}{\bar\eta}\arrow{s}\node[2]{\mathrm{Spin}(V,g)\subset \mathrm{Spin}(M,g)}\arrow{s}\\   
\node{U\subset T_pM}\arrow[2]{e,t}{\exp_p}\node[2]{V\subset M}   
\end{diagram}$$   
Hence, we obtain a map between the spinor bundles $\Sigma U$ and $\Sigma V$ in the following way:   
\begin{eqnarray}\label{indentspinbun}   
\Sigma U=\mathrm{Spin}(U,\geucl)\times_\rho\Sigma_n &\longrightarrow& \Sigma V=\mathrm{Spin}(V,g)\times_\rho\Sigma_n\nonumber\\   
\psi=[s,\varphi]&\longmapsto&\overline{\psi}=[\bar\eta(s),\varphi]   
\end{eqnarray} 
where $(\rho,\Sigma_n)$ is the complex spinor representation, and where
$[s,\phi]$ denotes the equivalence class of $(s,\phi)$ under the diagonal
action of $\mathrm{Spin}(n)$.
   
We now define  
  $$e_i:= b_i^j\partial_j\;,$$  
so that $(e_1,\ldots,e_n)$ is an   
orthonormal frame of $(TV,g)$.  
Denote by $\nabla$ (resp. $\bar\nabla$) the Levi-Civita connection  
on $(TU,\geucl)$ (resp. $(TM,g)$) as well as its lift    
to the spinor bundle $\Sigma U$ (resp. $\Sigma V$).  
The Christoffel symbols of the second kind  $\widetilde\Gamma^k_{ij}$ 
are defined by    
  $$\widetilde{\Gamma}^k_{ij}:=\langle\bar\nabla_{e_i}e_j,e_k\rangle\:,$$ 
hence $\widetilde\Gamma^k_{ij}=-\widetilde\Gamma^j_{ik}$.   
   
\begin{prop}\label{relatdirac}   
If $D$ and $\bar D$ denote the Dirac operators acting respectively on   
$\Gamma(\Sigma U)$ and $\Gamma(\Sigma V)$, then we have   
\begin{equation}   
\bar
D\bar\psi=\overline{D\psi}+\mathbf{W}\cdot\bar\psi+\mathbf{V}\cdot\bar\psi
+ \sum_{ij} (b_i^j-\de_i^j) \overline{ \partial_i \cdot \nabla_{\partial j}
 \psi} \;,    
\end{equation}   
where $\mathbf{W}\in\Gamma(\mathrm{Cl\ } TV)$  
and $\mathbf{V}\in\Gamma(TV)$ are defined by   
\begin{equation}\label{eq.w}
\mathbf{W}=\frac{1}{4}\sum_{\ijkdiff}b_i^r(\partial_r
b_j^l)(b^{-1})^k_le_i\cdot e_j\cdot e_k
\end{equation}
and   
\begin{equation}\label{eq.v}
\mathbf{V}=\frac{1}{4}\sum_{i,k}\Big(\widetilde\Gamma^i_{ik}-\widetilde\Gamma^k_{ii}\Big)\,e_k=\frac{1}{2}\sum_{i,k}\widetilde\Gamma^i_{ik}\,e_k
\end{equation}   
where, for any point $m\in V$, and 
the coefficients of the inverse matrix of $B_m$ are denoted by 
$(b^{-1})^k_l(m)$ .   
\end{prop}   
  
\begin{proof}
We denote Clifford multiplication on $\Sigma V$ by ``$\;\cdot\;$''. For all    
spinor field $\psi\in\Gamma(\Sigma U)$, since   
$\bar\psi\in\Gamma(\Sigma V)$ and by definition of $\bar\nabla$
(see e.g.\ \cite[Theorem 4.14]{lawson.michelsohn:89}, 
\cite[I Lemma~4.1]{baer:thesis}), we have   
\begin{equation} 
\label{relatnabla} 
\bar\nabla_{e_i}\bar\psi=\overline{\nabla_{e_i}(\psi)}+\frac{1}{4}\sum_{j,k}   
\widetilde\Gamma_{ij}^k\,e_j\cdot   
e_k\cdot\bar\psi\;. 
\end{equation}   
Taking Clifford multiplication by $e_i$ on each member of \eqref{relatnabla} and   
summing over $i$ yields
\begin{equation*}   
\bar D\bar\psi=\sum_i e_i \cdot \overline{\nabla_{e_i} \psi}+
\frac{1}{4}\sum_{i,j,k}\widetilde\Gamma^k_{ij}e_i\cdot e_j\cdot e_k\cdot\bar\psi\;.    
\end{equation*}   
Now, using that $e_i = \sum_j b_i^j \pa_j $ and that 
$   e_i \cdot \overline{\nabla_{e_i} \psi} = \overline{\pa_i \cdot
  \nabla_{e_i} \psi}  $, we obtain that 
$$ \bar D\bar\psi=\sum_{ij} b_i^j \overline{\pa_i \cdot\nabla_{\pa_j} \psi}+
\frac{1}{4}\sum_{i,j,k}\widetilde\Gamma^k_{ij}e_i\cdot e_j\cdot
e_k\cdot\bar\psi\,$$ 
and hence,
$$ \bar D\bar\psi=\overline{D\psi}  + \sum_{ij} (b_i^j - \de_i^j) 
\overline{\pa_i \cdot  \nabla_{\pa_j} \psi}+
\frac{1}{4}\sum_{i,j,k}\widetilde\Gamma^k_{ij}e_i\cdot e_j\cdot
e_k\cdot\bar\psi.$$
See also \cite{pfaeffle:thesis} for a similar formula, 
worked out in more detail.  

Note  that by the definition of $e_k$, we have    
  $$\widetilde\Gamma^k_{ij}e_k=\widetilde\Gamma^k_{ij}b_k^l\partial_l\;.$$   
On the other hand, we compute the Christoffel symbols of the second kind     
  $$\widetilde\Gamma^k_{ij}e_k=\bar\nabla_{e_i}e_j=b_i^r\bar\nabla_{\partial_r}(b_j^s\partial_s) 
    =b_i^r(\partial_r b_j^s)\partial_s+b_i^r b_j^s    
\Gamma_{rs}^l\partial_l\;,$$   
where as usually the Christoffel symbols of the first kind $\Gamma_{rs}^l$  
are defined by $$\Gamma_{rs}^l\partial_l=\bar\nabla_{\partial_r}\partial_s\;.$$   
   
Therefore we have   
$$\widetilde\Gamma^k_{ij}b_k^l = b_i^r(\partial_r b_j^l)+b_i^r b_j^s    
\Gamma_{rs}^l\;,$$   
and hence   
\begin{equation}\label{gamma1}   
\widetilde\Gamma^k_{ij}=\Big(b_i^r(\partial_r b_j^l)+b_i^r b_j^s\Gamma_{rs}^l\Big)(b^{-1})^k_l\;.   
\end{equation}   
   
Now, we can write   
$$\frac{1}{4}\sum_{i,j,k}\widetilde\Gamma^k_{ij}e_i\cdot e_j\cdot e_k=\mathbf{W}+\mathbf{V}$$   
where $\mathbf{W}\in\Gamma(\Lambda^3 TV)$ and $\mathbf{V}\in\Gamma(TV)$ are defined by   
$$\mathbf{W}=\frac{1}{4}\sum_{\ijkdiff}\widetilde\Gamma^k_{ij}e_i\cdot e_j\cdot e_k$$   
and   
\begin{eqnarray*}   
\mathbf{V}&=&\frac{1}{4}\Big(\sum_{i=j\not=k}\widetilde\Gamma^k_{ij}e_i\cdot
e_j\cdot e_k+\overbrace{\sum_{i\not=j=k}\widetilde\Gamma^k_{ij}e_i\cdot
  e_j\cdot e_k}^{=0}+\sum_{j\not=i=k}\widetilde\Gamma^k_{ij}e_i\cdot e_j\cdot
e_k+
\overbrace{\sum_{i=j=k}\widetilde\Gamma^k_{ij}e_i\cdot e_j\cdot e_k}^{=0}\Big)\\        
&=&\frac{1}{4}\sum_{i,k}\Big(\widetilde\Gamma^i_{ik}-\widetilde\Gamma^k_{ii}\Big)e_k   
\end{eqnarray*}  
which is \eref{eq.v}. 
   

First note that by \eqref{gamma1} we have    
$$\mathbf{W}=\frac{1}{4}\sum_{\ijkdiff}\Big(b_i^r(\partial_r b_j^l)(b^{-1})^k_l+b_i^r b_j^s\Gamma_{rs}^l(b^{-1})^k_l\Big)\;e_i\cdot e_j\cdot e_k\;.$$   
However,    
  $$\sum_{\ijkdiff}b_i^r b_j^s\Gamma_{rs}^l(b^{-1})^k_l\;e\cdot   
    e_j\cdot e_k=0$$    
since  $\Gamma_{rs}^l=\Gamma_{sr}^l$ and $e_i\cdot e_j=- e_j\cdot e_i.$   
Therefore  we obtain \eref{eq.w}.
\end{proof}


\section{Development of the metric at the point $p$} \label{dev}   
   
In this section we give the development of the coefficients $\widetilde\Gamma^k_{ij}$ in the coordinates $(x_1,\ldots,x_n)$ at    
the fixed point $p\in M$.    
   
For any point $m \in M$, $r$ denotes the distance from $p$ to $m$. Recall that in the neighborhood of $p$, we have the following    
development of the metric $g$ (see for example \cite{lee.parker:87}):   
   
\begin{eqnarray}\label{devmetric}   
g_{ij}&=& \delta_{ij}+    
\frac{1}{3} R_{i \alpha \beta j}(p) x^{\alpha} x^{\beta}   
+\frac{1}{6} R_{i \alpha \beta j;\gamma}(p) x^{\alpha} x^{\beta} x^{\gamma}\\   
&&+\left( \frac{1}{20}  R_{i \alpha \beta j;\gamma\lambda} (p)    
+ \frac{2}{45} \sum_m  R_{i \alpha \beta m}(p) R_{j \gamma \lambda   
  m}(p)\right) x^{\alpha} x^{\beta} x^{\gamma} x^{\lambda} + O(r^5)   
 \nonumber   
\end{eqnarray}   
   
where    
  $$R_{i j k l}= \langle \na_{e_j}\na_{e_i} e_k, e_l\rangle  
    - \langle \na_{e_i}\na_{e_j} e_k, e_l\rangle  
    - \langle \na_{[e_j,e_i]} e_k, e_l\rangle$$ 
and where  
  $$R_{i j k l;m}= (\nabla R)_{m i j k l }\qquad\qquad R_{i j k l;m n}=
  (\nabla^2R)_{n m i j k l }$$   
are the covariant derivatives of $R_{ijkl}$ in direction of $e_m$ (and $e_p$). 
Therefore we write $$G_m=\mathrm{Id}+G_2+G_3+O(r^4)$$ with   
$$\Big(G_2\Big)_{ij}=\frac{1}{3}R_{i \alpha \beta j}(p) x^{\alpha} x^{\beta}$$   
and   
$$\Big(G_3\Big)_{ij}=\frac{1}{6} R_{i \alpha \beta j;\gamma}(p) x^{\alpha} x^{\beta} x^{\gamma}$$   
   
\noindent Writing $$B_m=\mathrm{Id}+B_1+B_2+B_3+O(r^4)$$   
with $$\Big(B_1\Big)_{ij}=B_{ij\alpha}x^\alpha\;,$$   
$$\Big(B_2\Big)_{ij}=B_{ij\alpha\beta}x^\alpha x^\beta$$   
and   
$$\Big(B_3\Big)_{ij}=B_{ij\alpha\beta\gamma}x^\alpha x^\beta x^\gamma\;,$$   
the relation $B_m^2G_m= \mathrm{Id}$ yields   $B_1=0$ and 
$$0=\Big(2B_2+G_2\Big)+\Big(2B_3+G_3\Big)\;,$$   
hence   
\begin{equation}\label{devbij}   
\left\{\begin{array}{rcc}   
b_i^j&=&\delta_i^j-\frac{1}{6}R_{i\alpha\beta j}x^\alpha x^\beta-\frac{1}{12}R_{i\alpha\beta j;\gamma}x^\alpha x^\beta    
x^\gamma+O(r^4)\\[3mm]   
(b^{-1})_i^j&=&\delta_i^j+\frac{1}{6}R_{i\alpha\beta j}x^\alpha x^\beta+\frac{1}{12}R_{i\alpha\beta j;\gamma}x^\alpha x^\beta    
x^\gamma+O(r^4)   
\end{array}\right.   
\end{equation}   
We also have   
\begin{equation}\label{devbijder}   
\partial_lb_i^j=-\frac{1}{6}\Big(R_{il\alpha j}+R_{i\alpha lj}\Big)x^\alpha-\frac{1}{12}\Big(R_{il\alpha j;\beta}+R_{i\alpha lj;\beta}+R_{i\alpha\beta j;l}\Big) x^\alpha x^\beta+O(r^3)\;.   
\end{equation}   
   
\medskip   
   
\subsection{Development of $\Gamma_{ij}^k$, $\mathbf{V}$ and $\mathbf{W}$}   
   
\begin{eqnarray*}   
 \Gamma_{ij}^k&=&\frac{1}{2}g^{kl}\Big(\partial_ig_{jl}+\partial_jg_{il}-   
 \partial_lg_{ij}\Big)\\   
&=&\frac{1}{2}\Big(\partial_ig_{jk}+\partial_jg_{ik}-\partial_kg_{ij}\Big)+O(r^2)\\   
&=&\frac{1}{6}\Big(R_{ji\alpha k}+R_{j\alpha ik}+R_{ij\alpha k}+R_{i\alpha   
  jk}-R_{ik\alpha j}-R_{i\alpha kj}\Big)x^\alpha+O(r^2)   
\end{eqnarray*}   
   
Using the relations      
$ R_{ij\alpha k}+ R_{ji\alpha k}=0$, $R_{j\alpha ik}-R_{ik\alpha
  j}=-2R_{ik\alpha j}$ and $R_{i\alpha jk}-R_{i\alpha kj}=-2R_{i\alpha kj}$
we then have    
   
\begin{eqnarray}\label{devgammaijk}   
\Gamma_{ij}^k&=&-\frac{1}{3}\Big(R_{ik\alpha j}+R_{i\alpha   
  kj}\Big)x^\alpha+O(r^2)   
\end{eqnarray}   
   
On the other hand, since $\partial_rb_j^l$ and $\Gamma^l_{rs}$ 
have no   
constant term, Formula \eqref{gamma1} yields   
   
$$\widetilde\Gamma^k_{ij}=\Big(\delta_i^r(\partial_r b_j^l)+\delta_i^r   
\delta_j^s\Gamma_{rs}^l\Big)\delta^k_l+O(r^2)\;,$$   
   
and hence   
   
$$\widetilde\Gamma^k_{ij}=\partial_i b_j^k+\Gamma^k_{ij}+O(r^2)\;.$$

\noindent We have   
\begin{eqnarray*}   
  \mathbf{V}&=&\frac{1}{4}\sum_{i,k}\Big(\widetilde\Gamma^i_{ik}-\widetilde\Gamma^k_{ii})e_k\\   
&=&\frac{1}{4}\sum_{i,k}\Big(\partial_i b_k^i+\Gamma^i_{ik}-\partial_i   
b_i^k-\Gamma^k_{ii}\Big)e_k\\   
&=&\frac{1}{4}\sum_{i,k}\Big(\Gamma^i_{ik}-\Gamma^k_{ii}\Big)e_k   
\end{eqnarray*}   
since $\partial_i b_k^i=\partial_i b_i^k$.

\noindent Moreover, we have   
\begin{eqnarray*}   
 \sum_i\Big(\Gamma^i_{ik}-\Gamma^k_{ii}\Big) &=& 
-\frac{1}{3}\sum_i\Big(R_{ii\alpha k}+R_{i\alpha   
  ik}\Big)x^\alpha 
+\frac{1}{3}\sum_i\Big(R_{ik\alpha i}+R_{i\alpha   
  ki}\Big)x^\alpha+ O(r^2)\\   
&&=-(\Ric)_{\alpha k}+O(r^2)   
\end{eqnarray*} 
   
\noindent Therefore we proved that   
\begin{equation}   
  \label{devV}   
 \boxed{ \mathbf{V}=\Big(-\frac{1}{4}(\Ric)_{\alpha k}\,x^\alpha+   
   O(r^2)\Big)e_k\;.}    
\end{equation}   
   

The aim now is to show that   
 $$\mathbf{W}=\frac{1}{4}\sum_{\ijkdiff}b_i^r(\partial_r   
b_j^l)(b^{-1})^k_le_i\cdot e_j\cdot e_k\;$$ 
is $O(r^3)$. First note that by Equations (\ref{devbij}) and  
(\ref{devbijder}) $b_i^r$ has no term of order $1$ and 
$\partial_r b_j^l$ has no term of order $0$.  
Hence, any term in $\mathbf{W}$ of order $<3$ is a product of the $0$-order  
term of $b_i^r$ and of a term of order $1$ or $2$ of $\partial_r b_j^l$. 
   
   
Therefore $\mathbf{W}$ has no term of order $0$. To compute the terms of   
order $1$ and $2$, we write  
  $$\mathbf{W}=\frac{1}{4}\sum_{\ijkdiff}\Big(\delta_i^r(\partial_r   
b_j^l)\delta^k_l+O(r^3)\Big)e_i\cdot e_j\cdot e_k\;.$$  
We have   
  $$\sum_{\ijkdiff}\partial_ib_j^k e_i\cdot e_j\cdot e_k=0$$   
since   
  $$\partial_ib_j^k =\partial_ib_k^j\qquad\text{and}\qquad e_j\cdot   
e_k=-e_k\cdot e_j\;.$$   
   
Therefore $\mathbf{W}$ has no term of order $1$ and $2$. We proved that   
\begin{eqnarray}    
 \label{devW}   
\mathbf{W} =O(r^3)   
\end{eqnarray}   

\begin{rem}\label{rem.W}   
Similar calculations yield  
\begin{equation*}   
 \boxed{ \mathbf{V}=-\Big(\frac{1}{4}(\Ric)_{\alpha k}\,x^\alpha+\frac{1}{6}(\Ric)_{\alpha k,\beta}\,x^\alpha   
x^\beta + O(r^3)\Big)e_k\;.}   
\end{equation*}   
   
\begin{equation*}   
\boxed{ \mathbf{W}= -\frac{1}{144}\sum_{\ijkdiff}R_{l\beta\gamma   
  k}\Big(R_{ji\alpha l}+R_{jl\alpha i}\Big)\,x^\alpha x^\beta x^\gamma\, e_i\cdot e_j\cdot   
e_k+O(r^4)\;.}   
\end{equation*}   
We do not give details here because we do not need explicit computations of   
terms of order 2 for $\mathbf{V}$ and terms of order 3 for $\mathbf{W}$ in   
the proof of Theorem~\ref{main}.   
\end{rem}   
   
\section{The test spinor \label{testspinor}}

\subsection{The explicit spinor}\label{subsec.explicit}
In this section we construct a good test spinor on $\mR^n$. The spinor bundle on $\mR^n$ is trivial, 
so we can identify the fibers. Let $\psi_0\in \Sigma_0\mR^n$.
We set $f(x):=\frac{2}{1+r^2}$, where $r:=|x|$, hence $\pa_i f= -x_i f^2$.  
Then we define
\begin{eqnarray}\label{explicit}
\psi(x) =  f^{\frac{n}{2}}(x)(1 - x) \cdot \psi_0.
\end{eqnarray}
One calculates
\begin{equation}\label{few0} 
\nabla_{\pa_i} \psi = -f^{\frac{n}{2}}\pa_i\cdot \psi_0 - \frac{n}2 f^{\frac{n}{2}+1}x_i(1-x)\cdot\psi_0,
\end{equation}
and hence
\begin{eqnarray}   
  D\psi&=&\frac{n}{2}f\psi\label{few}\\   
|\psi|&=&f^\frac{n-1}{2}\label{few2}\\   
|D\psi|&=&\frac{n}{2}f^\frac{n+1}{2}.\label{few3}   
\end{eqnarray}


\subsection{Conformal change of metrics}\label{subsec.confch}

In order to explain a geometric interpretation of 
this spinor, we have to recall the behavior 
of spinors and the Dirac operators under conformal changes. See e.g.\
\cite{hitchin:74,hijazi:01} for proofs.

Let $(N,h)$ be a spin manifold of dimension $n$.   
Consider a conformal change of metric $\widetilde{h}=F^{-2}h$ for any 
positive real function $F$ on $(N,h)$. The map $TN\to TN$, 
$X\mapsto \tilde X= FX$ induces
an isomorphism of principal bundles from $\mathrm{SO} (N,h)$ to $\mathrm{SO} 
(N,\widetilde{h})$. It lifts to a bundle isomorphism 
between the $\mathrm{Spin}(n)$-principal bundles    
$\mathrm{Spin} (N,h)$ and $\mathrm{Spin} (N,\widetilde{h})$, and 
passing to the associated bundles one obtains a map 
\begin{eqnarray*}
 \Sigma_h N=  \mathrm{Spin} (N,h)\times_\rho \Sigma&\to&
\Sigma_{\widetilde h} N=  \mathrm{Spin} (N,\widetilde h)\times_\rho \Sigma\\
\phi& \mapsto & \widetilde{\phi}
\end{eqnarray*}
between the spinor bundles, which is a fiberwise isometry and we have
 $$\widetilde{X}\,\tilde\cdot\,\widetilde{\phi}=\widetilde{X\cdot\phi}$$     
(see \cite{hijazi:01} for more details on this construction). 

By conformal covariance of the Dirac operator, we have, for $\varphi \in    
 \Gamma(\Sigma N)$,   
\begin{equation}\label{1covconf}   
\widetilde{D}\Big(\; F^{\frac{n-1}{2}} \;\widetilde{\varphi}\; \Big)=F^{\frac{n+1}{2}} \;\widetilde{D \varphi},     
\end{equation}   
   
 
\subsection{Geometric interpretation}

We apply this formula to a particular case: let $p$ be any point of the round sphere $\mS^n$. Then $\mS^n\backslash \{p\}$ is isometric to $\mR^n$    
with the metric  
\begin{eqnarray} \label{conform}  
g_S=f^2g_{\text{eucl}}\;, 
\end{eqnarray}  
 with $$f(x)=\frac{2}{1+r^2}\;.$$   
Hence we set $N:=\mR^n$, $h=\geucl$, $F=f^{-1}$.
One calculates with \eref{few} and \eref{1covconf}
that $\Phi:=F^{\frac{n-1}2}\tilde \psi$ satisfies $D\Phi=\frac{n}2 \Phi$ on $\mS^n\setminus\{p\}$, and $|\Phi|=1$.
Hence, the possible singularity at $p$ can be removed (see e.g. the Removal of singularity theorem
\cite[Theorem~5.1]{ammann:p03}), and one sees that $\Phi$ is 
an eigenspinor to the eigenvalue $n/2$ on the round sphere $\mS^n$. 
The equality discussion in Friedrich's inequality 
\cite{friedrich:80} implies that $\Phi$ is a Killing 
spinor to the constant $-1/2$,
i.e.\ it satisfies
  $$\na_X\Phi =-\frac12 X\cdot \Phi.$$
Hence we have seen that our spinor $\psi$  is the ``conformal image'' of a Killing spinor on $\mS^n$.

\section{The proof of Theorem \ref{main} for $n\geq 3$} \label{estimate}   

We begin with the following Proposition.   
   
\begin{prop}\label{choixmetric} The metric $g$ on $M$ can be chosen such that   
$$\Ric_g(p)=0\qquad\text{and}\qquad \Delta_g(\Scal_g)(p)=0\;.$$   
\end{prop}   
   
\begin{proof}   
Consider a conformal change of the metric $\widetilde{g}=e^{2u}g$ for any real    
function $u$ on $(M,g)$. Then it is well known that the Ricci curvature $(2,0)$-tensor   
$\Ric_{\tilde{g}}$, the scalar curvature   
$\Scal_{\tilde{g}}$ and the Laplacian $\Delta_{\tilde{g}}$   
corresponding to the metric ${\tilde{g}}$ satisfy (see for example Hebey   
\cite{hebey:97} or Aubin \cite{aubin:76})   
$$\Ric_{\tilde{g}}=\Ric_g-(n-2)\nabla^2u+(n-2)\nabla   
u\otimes\nabla u+(\Delta_g u-(n-2)|\nabla u|^2_g)g\;,$$   
\begin{equation}   
  \label{relatscal}   
\Scal_{\tilde{g}}=e^{-2u}\Big(\Scal_g+2(n-1)\Delta_gu-(n-1)(n-2)|\nabla   
u|^2_g\Big)\;,     
\end{equation}   
   
\noindent As a first step, we can assume that  
$\Scal_g(p)=0$. Then, let us choose $u$ such  that $$u(x)=   
\frac{1}{2(n-2)}\Big(\Ric_g(p)_{ij}-\frac{\Scal_g(p)}{n}g_{ij}(p)\Big)\,x^ix^j-\frac{\Delta_g(\Scal_g)(p)}{48(n-1)}(x^1)^4$$   
in a neighborhood of the point $p$. Since $u(p)=0$ and $(\nabla u)(p)=0$,   
it is straightforward to see that $\Ric_{\tilde g}(p)=0$. Moreover,   
taking the Laplacian of both members of Equation \eqref{relatscal}, a   
simple computation shows that   
$\Delta_{\tilde{g}}\Scal_{\tilde{g}}(p)=0$.    
\end{proof}   
   
Let $\bar\varphi\in\Sigma_UM$ where $U$ is the open neighborhood of a point   
$p\in M$ as defined in the previous sections. With the help of formulas   
\eqref{devV} and \eqref{devW}, we have the following    
   
\begin{cor} \label{coro}
For any metric $g$ on $M$ chosen as in Proposition \ref{choixmetric},   
 we have    
\begin{eqnarray}   
\label{relatDDbar}   
\bar   
D\bar\varphi=\overline{D\varphi}&+&\sum_{\ijkabc}A_{ijk\alpha\beta\gamma}   
x^\alpha\,x^\beta\,x^\gamma\, e_i\cdot e_j\cdot e_k \cdot\bar\varphi +   
\mathbf{W}'\cdot \bar\varphi+\mathbf{V} \cdot \bar\varphi+ \sum_{ij} (b_i^j-\de_i^j) \overline{ \partial_i \cdot \nabla_{\partial j}
 \psi} \;\\   
 \end{eqnarray}   
where $A_{ijk\alpha\beta\gamma}\in\mR$
and where   
$\mathbf{W}' \in\Gamma(\Lambda^3 TV)$, $\mathbf{V}\in\Gamma(TV)$,   
$|\mathbf{W}'|\leq C\, r^4$ and $|\mathbf{V} |\leq C'\, r^2$    
($C$ and $C'$ being positive constants independent of $\varphi$).    
\end{cor}   

\begin{rem} \label{coro.remark} Using the formulae in Remark~\ref{rem.W}, 
we obtain 
the formula
  $$A_{ijk\alpha\beta\gamma}=-\frac{1}{144}R_{l\beta\gamma   
  k}\Big(R_{ji\alpha l}+R_{jl\alpha i}\Big)\;,$$    
\end{rem}   
    
Assume now that $\psi$ is the test spinor constructed in   
Section \ref{testspinor}. Let $\ep >0$ be a small positive number.   
We set    
$$\varphi(x):=\eta\psi(\frac{x}{\ep})=:\psi_\ep(x)$$ where $\eta=0$   
on $\mR^n\setminus B_p(2\delta)$ and $\eta=1$ on $B_p(\delta)$, and that $\psi$, defined as in \eref{explicit}  
satisfies the following relations \eref{few0}, \eref{few}, \eref{few2} and \eref{few3}
where $f$ is again defined by $$f(x)=\frac{2}{1+r^2}\;.$$

We now prove some lemmas which will be useful in the proof of Theorem
\ref{main}. 
   
\begin{lem} \label{term_supp}
We have
\begin{eqnarray} \label{es1}
  \left| \sum_{ij}  (b_i^j- \de_i^j) \pa_i \cdot
    \nabla_{\pa_j}(\psi(\frac{x}{\ep}))  \right| \leq C \frac{r^3}{\ep} f^{\frac{n}{2}}(\frac{x}{\ep})
\end{eqnarray} 
where $f=\frac{2}{1+r^2};$.
\end{lem}    

\begin{proof}

At first, we prove that: 
\begin{eqnarray} \label{devo3} 
\sum_{ij\al\be} R_{i\al \be j} x^\al x^\be \pa_i \cdot   \left( \nabla_{\pa_j}
\psi(\frac{x}{\ep}) \right) =0.
\end{eqnarray}

Indeed, using (\ref{few0}), we compute that 
$$  (\nabla_{\pa_j}\psi)(\frac{x}{\ep})=  - \frac{f^{\frac{n}{2}} (\frac{x}{\ep})}{\ep}
\pa_j \cdot \psi_0 -  \frac{ n f^{\frac{n+2}{2}} (\frac{x}{\ep})}{2\ep} 
 x^j (1-\frac{x}{\ep}) \cdot  \psi_0.$$
and obtain
$$\sum_{ij\al\be} R_{i\al \be j} x^\al x^\be \pa_i \cdot   (\nabla_{\pa_j}
\psi)(\frac{x}{\ep})=  - \frac{f^{\frac{n}{2}} (\frac{x}{\ep})}{\ep}
\sum_{ij\al\be} R_{i\al \be j} x^\al x^\be  
\pa_i \cdot \pa_j \cdot \psi_0 -  \frac{ n
  f^{\frac{n+2}{2}} (\frac{x}{\ep})}{2\ep} 
\sum_{ij\al\be} R_{i \al \be j}  x^\al x^\be x^j \pa_i \cdot (1-\frac{x}{\ep})
\cdot  \psi_0.$$

Now, since if $i \not= j$, $\pa_i \cdot \pa_j= -\pa_j \cdot \pa_i$ and
since $$\sum_{\al \be}  R_{i\al \be j}  x^\al x^\be = 
\sum_{\al \be} R_{i \be \al j} x^\al x^\be = 
\sum_{\al \be} R_{j\al  \be i} x^\al x^\be,$$
(we have used that $R_{j \al \be i}= R_{\be i j \al} = R_{i \be \al j}$),
we get that 
$$\sum_{ij\al\be} R_{i\al \be j} x^\al x^\be 
\pa_i \cdot \pa_j \cdot \psi_0 = - \sum_{i,
  \al \be}  R_{i\al \be i} x^\al x^\be \psi_0 = 0$$
since $\Ric(p)=0$.
The first summand vanishes.

The second summand vanishes as $\sum_{\beta j}R_{i\al\be j}x^\be x^j=0$.

This proves (\ref{devo3}). Now, by the development of $b_i^j$
(\ref{devbij}), we easily obtain that 
$$  \left| \sum_{ij}  (b_i^j- \de_i^j) \pa_i \cdot
    \nabla_{\pa_j}(\psi(\frac{x}{\ep}))  \right| \leq C \frac{r^3}{\ep}
|\nabla \psi|(\frac{x}{\ep}).$$ 

Differentiating expression (\ref{explicit}), we see that 
$$|\nabla \psi| \leq C( f^{\frac{n}{2}} +  r f^{\frac{n+2}{2}}).$$
Together with $rf(r) = \frac{2r}{1+r^2}= 1 - \frac{(1-r)^2}{1+r^2}\leq 1$
we obtain the lemma.
\end{proof}

Now, we can start the proof of Theorem~\ref{main}. We have, with
the notations of Corollary~\ref{coro}:    
\begin{eqnarray*}   
  \bar   D{\bar\psi_{\ep}(x)}&=&\bar\nabla\eta\cdot{\bar\psi(\frac{x}{\ep})}+\eta\,\bar D({\bar\psi(\frac{x}{\ep})})\\&=&\bar\nabla\eta\cdot\bar\psi(\frac{x}{\ep})+\frac{\eta}{\ep}\,\overline{ D\psi}(\frac{x}{\ep})+\eta\,\sum_{\ijkabc}A_{ijk\alpha\beta\gamma}   
x^\alpha\,x^\beta\,x^\gamma\, e_i\cdot e_j\cdot e_k
\cdot\bar\psi(\frac{x}{\ep})\\   
&&+\eta\,\mathbf{W}' \cdot \bar\psi(\frac{x}{\ep})+   
\eta\mathbf{V} \cdot \bar\psi(\frac{x}{\ep})+ \eta\sum_{ij} (b_i^j-\de_i^j) 
\overline{ \partial_i \cdot \nabla_{\partial j}(\psi(\frac{x}{\ep}))}\;.   
\end{eqnarray*} 

Therefore we have   
\begin{eqnarray*}   
 \bar D\bar\psi_{\ep}(x)&=&\bar\nabla\eta\cdot\bar\psi(\frac{x}{\ep})+\frac{\eta}{\ep}\,\frac{n}{2}\, f(\frac{x}{\ep})\,\bar\psi(\frac{x}{\ep})+\eta\,\sum_{\ijkabc}A_{ijk\alpha\beta\gamma}   
x^\alpha\,x^\beta\,x^\gamma\, e_i\cdot e_j\cdot e_k   
\cdot\bar\psi(\frac{x}{\ep})\\   
&&+\eta\,\mathbf{W}' \cdot \bar\psi(\frac{x}{\ep})+   
\eta\mathbf{V} \cdot \bar\psi(\frac{x}{\ep}) +  \eta \sum_{ij} (b_i^j-\de_i^j) \overline{ \partial_i \cdot \nabla_{\partial j}
 (\psi(\frac{x}{\ep}))}\;.   
\end{eqnarray*}

We write that
\begin{eqnarray*} 
|\bar D\bar \psi_\ep|^2(x)& =&   
\mathbf{I}+\mathbf{II}+\mathbf{III}+\mathbf{IV}+\mathbf{V}+\mathbf{VI}+\mathbf{VII}+\mathbf{VIII}+\mathbf{IX}+\mathbf{X}+\mathbf{XI}+\mathbf{XII}+\mathbf{XIII}+\mathbf{XIV}+\mathbf{XV}\\
&&+\mathbf{XVI}+\mathbf{XVII}+\mathbf{XVIII}+\mathbf{XIX}+\mathbf{XX}+\mathbf{XXI}\;
\end{eqnarray*}
where    
\begin{eqnarray*}   
\mathbf{I}&=&|\bar\nabla\eta|^2\,|\bar\psi|^2(\frac{x}{\ep})\\
\mathbf{II}&=&\frac{\eta^2}{\ep^2}\,\frac{n^2}{4}\,|\bar\psi|^2(\frac{x}{\ep})\,f^2(\frac{x}{\ep})\\   
\mathbf{III}&=&\eta^2|\sum_{\ijkabc}A_{ijk\alpha\beta\gamma}   
x^\alpha\,x^\beta\,x^\gamma\, e_i\cdot e_j\cdot e_k \cdot\bar\psi|^2(\frac{x}{\ep})\\   
\mathbf{IV}&=&\eta^2\,|{\mathbf{W}'}|^2
|\bar\psi(\frac{x}{\ep})|^2\\  
 \mathbf{V} & =& \eta^2 |\mathbf{V}|^2 |\bar\psi(\frac{x}{\ep})|^2\\   
\mathbf{VI}&=&2\,\Re   
e<\bar\nabla\eta\cdot\bar\psi(\frac{x}{\ep}),\frac{\eta}{\ep}\,\frac{n}{2}\,f(\frac{x}{\ep})\bar\psi(\frac{x}{\ep})>\\   
\mathbf{VII}&=&2\,\Re   
e<\bar\nabla\eta\cdot\bar\psi(\frac{x}{\ep}),\eta\sum_{\ijkabc}A_{ijk\alpha\beta\gamma}   
x^\alpha\,x^\beta\,x^\gamma\, e_i\cdot e_j\cdot e_k   
\cdot\bar\psi(\frac{x}{\ep})>\\   
\mathbf{VIII}&=&2\,\Re   
e<\bar\nabla\eta\cdot\bar\psi(\frac{x}{\ep}),\eta\,\mathbf{W}'   
\cdot \bar\psi(\frac{x}{\ep})>\\   
\mathbf{IX}&=&2\,\Re   
e<\bar\nabla\eta\cdot\bar\psi(\frac{x}{\ep}),\eta\,\mathbf{V}   
\cdot \bar\psi(\frac{x}{\ep})>\\   
\mathbf{X}&=&\frac{\eta^2}{\ep}\, n\,f(\frac{x}{\ep})\eta\sum_{\ijkabc}A_{ijk\alpha\beta\gamma}   
x^\alpha\,x^\beta\,x^\gamma\,\Re e <e_i\cdot e_j\cdot e_k   
\cdot\bar\psi,\bar\psi>(\frac{x}{\ep})\\   
\mathbf{XI}&=&\frac{\eta^2}{\ep}\, n\,f(\frac{x}{\ep})\Re e   
<\bar\psi(\frac{x}{\ep}),\mathbf{W}'\cdot \bar\psi(\frac{x}{\ep})>\\   
\mathbf{XII}&=&\frac{\eta^2}{\ep}\, n\,f(\frac{x}{\ep})\Re e   
<\bar\psi(\frac{x}{\ep}),\mathbf{V}\cdot   
\bar\psi(\frac{x}{\ep})>\\   
\mathbf{XIII}&=&2\,\eta^2\sum_{\ijkabc}A_{ijk\alpha\beta\gamma}   
x^\alpha\,x^\beta\,x^\gamma\,\Re e <e_i\cdot e_j\cdot e_k   
\cdot\bar\psi(\frac{x}{\ep}),\mathbf{W}'\cdot   
\bar\psi(\frac{x}{\ep})>\\   
\mathbf{XIV}&=&2\,\eta^2\sum_{\ijkabc}A_{ijk\alpha\beta\gamma}   
x^\alpha\,x^\beta\,x^\gamma\,\Re e <e_i\cdot e_j\cdot e_k   
\cdot\bar\psi(\frac{x}{\ep}),\mathbf{V} \cdot   
\bar\psi(\frac{x}{\ep})>\\   
\mathbf{XV}&=& 2 \eta^2 \Re e < \mathbf{W}'\cdot   
\bar\psi(\frac{x}{\ep}),\mathbf{V}\cdot   
\bar\psi(\frac{x}{\ep})>   \\
\mathbf{XVI}&=&  2 \Re e <\bar \nabla
\eta\cdot\bar\psi(\frac{x}{\ep}),   \eta \sum_{ij} (b_i^j-\de_i^j) \overline{ \partial_i \cdot \nabla_{\partial j}
 (\psi(\frac{x}{\ep}))}> \\    
 \mathbf{XVII}&=&    \frac{n\eta^2}{\ep}  f(\frac{x}{\ep})
 \Re e < \bar\psi(\frac{x}{\ep}), \sum_{ij} (b_i^j-\de_i^j) \overline{ \partial_i \cdot \nabla_{\partial j}
 (\psi(\frac{x}{\ep}))}> \\   
 \mathbf{XVIII}&=& 2\eta^2 \Re e<\sum_{\ijkabc}A_{ijk\alpha\beta\gamma}   
x^\alpha\,x^\beta\,x^\gamma\, e_i\cdot e_j\cdot e_k   
\cdot\bar\psi(\frac{x}{\ep}),\sum_{ij} (b_i^j-\de_i^j) \overline{ \partial_i \cdot \nabla_{\partial j}
 (\psi(\frac{x}{\ep}))}>  \\
\mathbf{XIX}&=&  2 \eta^2 \,\Re e< \mathbf{W}' \cdot
\bar\psi(\frac{x}{\ep}), \sum_{ij} (b_i^j-\de_i^j) \overline{ \partial_i \cdot \nabla_{\partial j}
 (\psi(\frac{x}{\ep}))}>    \\  
 \mathbf{XX}&=& 2 \eta^2 \,\Re e< \mathbf{V} \cdot
\bar\psi(\frac{x}{\ep}), \sum_{ij} (b_i^j-\de_i^j) \overline{ \partial_i \cdot \nabla_{\partial j}
 (\psi(\frac{x}{\ep}))}>    \\  
\mathbf{XXI}&=& \eta^2\left| \sum_{ij} (b_i^j-\de_i^j) \overline{ \partial_i \cdot \nabla_{\partial j}
 (\psi(\frac{x}{\ep}))} \right|^2.
\end{eqnarray*}
\noindent Since $\mathbf{V}$ is a vector field, we have    
$$\mathbf{XII}=0$$   
Assume now that $x\in B_p(2\delta)$. Using the fact that    
$|\bar\nabla \eta| \leq C r^4$  
($C$ being a constant independent of   
$\ep$) and since $r\leq\delta\leq 1$, we have:

$$\mathbf{I}+\mathbf{III}+\mathbf{IV}+\mathbf{V}+\mathbf{VII}+\mathbf{VIII}+   
\mathbf{IX}+\mathbf{XIII}+\mathbf{XIV}+\mathbf{XV}   
\leq C\,r^4\,f^{n-1}(\frac{x}{\ep})\;.$$   
and   
$$\mathbf{VI}+\mathbf{XI} \leq {C \over \ep}\,r^4\,f^{n}(\frac{x}{\ep})\;.$$   
Since $f \leq 2$ and since $r^2 \leq C$ on $B_p(2 \de)$, we obtain that 
 $$\mathbf{VI}+\mathbf{XI} \leq C\frac{r^2}{\ep} \, f^{n-\frac{1}{2}} (\frac{x}{\ep})\;.$$   
In the same way, using relation (\ref{es1}) and the fact that for all
$\ep$, 
$\frac{r}{\ep} f( \frac{x}{\ep}) \leq 1$, we have also 
$$\mathbf{X}+\mathbf{XVI}+\mathbf{XVII}+\mathbf{XVIII}+\mathbf{XIX}+\mathbf{XX}
 +
\mathbf{XXI}
\leq  C\frac{r^2}{\ep} \, f^{n-\frac{1}{2}} (\frac{x}{\ep})\;.$$ 

\noindent Therefore we obtain that   
\begin{eqnarray*}   
 |\bar D\bar \psi_\ep|^2(x)&\leq&\frac{n^2}{4\ep^2}   
f^{n+1}(\frac{x}{\ep})+
C\,r^4\,f^{n-1}(\frac{x}{\ep})+\frac{C}{\ep}\,r^2\,f^{n-\frac{1}{2}}(\frac{x}{\ep})\\   
&\leq&\frac{n^2}{4\ep^2}f^{n+1}(\frac{x}{\ep})\left[1+\Delta\right]   
\end{eqnarray*}   
   
\noindent where    
\begin{eqnarray*}   
  \Delta= C\,\ep^2 r^4\,f^{-2}(\frac{x}{\ep})+C \ep\,r^2\,f^{-\frac{3}{2}}(\frac{x}{\ep}).
\end{eqnarray*}   
   
\noindent Since $|\bar D\bar \psi_\ep|^2\geq 0$ we have $\Delta\geq   
-1$. Moreover, if we define   
$$g(x)=1+\frac{n}{n+1}x-(1+x)^\frac{n}{n+1}\;,\qquad\forall x\geq   
-1\;,$$ then   
$$g'(x)=\frac{n}{n+1}\Big(1-(1+x)^\frac{-1}{n+1}\Big)\;,\qquad\forall x>   
-1\;.$$ Therefore $g$ admits a minimum at $0$ on the interval   
$[-1,+\infty[$. This yields that, $\forall x\geq -1$,   
$$(1+x)^\frac{n}{n+1}\leq 1+\frac{n}{n+1}x\;.$$ We then have   
$$|\bar D\bar   
\psi_\ep|^\frac{2n}{n+1}(x)\leq(\frac{n}{2\ep})^\frac{2n}{n+1}\,f^n(\frac{x}{\ep})\,\left[1+\Delta\right]^\frac{n}{n+1}\leq(\frac{n}{2\ep})^\frac{2n}{n+1}\,f^n(\frac{x}{\ep})+\frac{n}{n+1}\,(\frac{n}{2\ep})^\frac{2n}{n+1}\,f^n(\frac{x}{\ep})\,\Delta\;.   
$$   
Taking into account the definition  of $\Delta$ and integrating over $M$   
leads to   
\begin{equation}   
  \label{majABCDEF}   
\int_M |\bar D\bar\psi_\ep|^\frac{2n}{n+1}\mathrm{d}v_g\leq\ep^{\frac{-2n}{n+1}}\left[\mathbf{A}+\mathbf{B}+\mathbf{C}\right]\;,   
\end{equation}   
where   
\begin{eqnarray*}   
\mathbf{A}&=& \int_{B_p(2\delta)}\left(\frac{n}{2}\right)^\frac{2n}{n+1}f^n(\frac{x}{\ep})\,\mathrm{d}v_g\\   
 \mathbf{B}&=&  C\,\int_{B_p(2\delta)} \ep^2 r^4\,f^{n-2}(\frac{x}{\ep})\,\mathrm{d}v_g  \\   
\mathbf{C}&=&    C\,\int_{B_p(2\delta)}\ep \,r^2\,f^{n-\frac{3}{2}}(\frac{x}{\ep})\,\mathrm{d}v_g.   
\end{eqnarray*}   
Since the function $f$ is radially symmetric, we can compute $\mathbf{A}$ with the help   
of spherical coordinates:   
$$\mathbf{A}=\int_{B_p(2\delta)}\left(\frac{n}{2}\right)^\frac{2n}{n+1}f^n(\frac{x}{\ep})\,\omega_{n-1}\,G(r)\,r^{n-1}\mathrm{d}r\;,$$   
where $\omega_{n-1}$ stands for the volume of the unit sphere $\mS^{n-1}$ and   
  $$G(r)=\frac{1}{\omega_{n-1}}\int_{\mS^{n-1}}\sqrt{|g|_{rx}} 
    \;\mathrm{d}\sigma(x)\;\qquad\qquad |g|_{y}:=\det g_{ij}(y).$$   
   
\noindent From Proposition \ref{choixmetric},  Hebey \cite{hebey:97} or Lee-Parker \cite{lee.parker:87}, we know that   
$$G(r)\leq 1+O(r^4)\;.$$  
Therefore, we can estimate   
$\mathbf{A}$ in the following way:   
\begin{eqnarray*}   
  \mathbf{A}&\leq   
  &\left(\frac{n}{2}\right)^\frac{2n}{n+1}\omega_{n-1}\left[\int_0^{2\delta}f^n(\frac{x}{\ep})\,r^{n-1}\mathrm{d}r+ C \int_0^{2\delta}f^n(\frac{x}{\ep})\,r^{n+3}\mathrm{d}r\right]\\   
&\leq &\left(\frac{n}{2}\right)^\frac{2n}{n+1}\omega_{n-1}\,\ep^n\left[\int_0^{\frac{2\delta}{\ep}}\frac{2^n 
    r^{n-1}}{(1+r^2)^n}\mathrm{d}r+C\ep^4\,\int_0^{\frac{2\delta}{\ep}}\frac{r^{n+3}}{(1+r^2)^n}\mathrm{d}r \right].  
\end{eqnarray*}   
Since
$$\int_0^{\frac{2\delta}{\ep}}\frac{r^{n+3}}{(1+r^2)^n}\mathrm{d}r
\leq O\left(
  \int_1^{\frac{2\delta}{\ep}} r^{3-n} \mathrm{d}r 
\right)$$ 
we get that 
$$\mathbf{A} \leq  \left(\frac{n}{2}\right)^\frac{2n}{n+1}\omega_{n-1}\,\ep^n\left[\int_0^{\frac{2\delta}{\ep}}\frac{2^n 
    r^{n-1}}{(1+r^2)^n}\mathrm{d}r+ o(1) \right]  \, ,$$
and hence   
\begin{equation}   
  \label{majorA}   
 \mathbf{A}\leq  \left(\frac{n}{2}\right)^\frac{2n}{n+1}\omega_{n-1}\,\ep^n\left[\int_0^{\frac{2\delta}{\ep}}\frac{2^n 
    r^{n-1}}{(1+r^2)^n}\mathrm{d}r+ o(1) \right]  \, .  
\end{equation}

\noindent Let us show that 
\begin{eqnarray} \label{es2} 
\mathbf{B}=o(\ep^n). 
\end{eqnarray}

\noindent Since  $dv_g \leq C dx$, setting $y=
\frac{x}{\ep}$, we have  
 \begin{eqnarray*} 
 \int_{B_p(2\delta)} r^4 f^{n-2} (\frac{x}{\ep}) \, dv_g \leq 
& \leq &  C \ep^{n+4} \int_{B_p(\frac{2 \delta}{\ep})} r^4 f^{n-2}
((\frac{y}{\ep})  \, dy \\
& \leq & C \,  \ep^{n+4}
  \int_0^\frac{2\delta}{\ep} \frac{ r^{n+3}}{(1+r^2)^{n-2}} \, dr.\\
& \leq & C \,  \ep^{n+4} O\left(  \int_1^\frac{2\delta}{\ep} r^{7-n} \, dr \right).
\end{eqnarray*}
It is easy to check that relation (\ref{es2}) follows if $n \geq 3$.
In the same way, we can prove that $\mathbf{C}=o(\ep^n)$.

\noindent Together with Equation  (\ref{majorA}),   
 we can conclude that   
\begin{eqnarray*}   
\int_M |\bar D\bar\psi_\ep|^\frac{2n}{n+1}\mathrm{d}v_g&\leq&\ep^{{\frac{-2n}{n+1}}+n}\left[\left(\frac{n}{2}\right)^\frac{2n}{n+1}\omega_{n-1}\int_0^{+\infty}r^{n-1}f^n(r)\mathrm{d}r +o(1)    
\right]\;,   
\end{eqnarray*}   
which yields   
\begin{equation}\label{presquemaj}   
\left(\int_M |\bar   
  D\bar\psi_\ep|^\frac{2n}{n+1}\mathrm{d}v_g\right)^\frac{n+1}{n}\leq\ep^{n-1}\left[\left(\frac{n^2}{4}\right)^\frac{n}{n+1}\omega_{n-1}I\right]^\frac{n+1}{n}    
\left( 1+o(1) \right)\;,   
\end{equation}   
where $$I=\int_0^{+\infty}\frac{2^n r^{n-1}}{(1+r^2)^n}\mathrm{d}r\;.$$

   
   
   
\noindent We are now going to estimate $\left|\int_M \Re e<\bar   
D\bar\psi_\ep,\bar\psi_\ep>\mathrm{d}v_g\right|$. We start by   
computing   
   

$$\left|\int_M \Re e<\bar   
  D\bar\psi_\ep,\bar\psi_\ep>\mathrm{d}v_g\right|\geq   
\mathbf{A'}-\mathbf{B'}-\mathbf{C'}-\mathbf{D'}-\mathbf{E'} \;,$$   
where   
\begin{eqnarray*}   
\mathbf{A'}&=&\int_{B_p(\delta)}\frac{n}{2\ep}f^n(\frac{x}{\ep})\mathrm{d}v_g\\
\mathbf{B'} & =& \left|  \int_M \Re
e<\bar\nabla\eta\cdot\bar\psi(\frac{x}{\ep}),\eta\bar\psi(\frac{x}{\ep})>
\mathrm{d}v_g\right| \,, \\
\mathbf{C'}& =& \left| \int_M\eta^2\sum_{\ijkabc}A_{ijk\alpha\beta\gamma}   
x^\alpha\,x^\beta\,x^\gamma\,\Re e <e_i\cdot e_j\cdot   
e_k\cdot\bar\psi,\bar\psi>(\frac{x}{\ep})\,\mathrm{d}v_g \right| \,
,  \\
\mathbf{D'} &=&  \left|\int_M\eta^2\Re
  e<\mathbf{W}'\bar\psi(\frac{x}{\ep}),\bar\psi(\frac{x}{\ep})>\,\mathrm{d}v_g \right| \\
\mathbf{E'} & =& \left|
\int_M  \eta^2   
 \Re e <  \sum_{ij} (b_i^j-\de_i^j) \overline{ \partial_i \cdot \nabla_{\partial j}
 (\psi(\frac{x}{\ep}))},\bar\psi(\frac{x}{\ep})
>\,\mathrm{d}v_g \right|  \, .
\end{eqnarray*}  
   
\noindent  (The term in $\mathbf{V}$ is zero). Note that    
$\mathbf{A'} =  \frac{1}{\ep} {\left( \frac{n}{2}\right)}^{1-   
  \frac{2n}{n+1}} \mathbf{A}$ where $\eta$ has been replaced by $2   
\eta$. 
As to obtain (\ref{es2}), we get that 
$$\mathbf{B'}+\mathbf{C'}+\mathbf{D'} \leq  C \int_{B_p(2\de)}
f^{n-1}(\frac{x}{\ep}) \,  \mathrm{d}v_g  \leq 0(\ep^{n}) = o(\ep^{-1})$$
and 
$$\mathbf{E'} \leq C \int_{B_p(2\de)}
\frac{r^3}{\ep} f^{n-\frac{1}{2}}(\frac{x}{\ep}) \,  \mathrm{d}v_g \leq o(\ep^{n-1}).$$  

Moreover, with the same method which was used to obtain (\ref{majorA}), we get

$$\mathbf{A'} \geq
\frac{n}{2}\,\omega_{n-1}\,\ep^{n-1}\,I\,\left[1+o(1) \right].$$

\noindent This proves that    
\begin{equation}   
  \label{mindpsipsi}   
\left|\int_M \Re e<\bar   
  D\bar\psi_\ep,\bar\psi_\ep>\mathrm{d}v_g\right|\geq
\frac{n}{2}\,\omega_{n-1}\,\ep^{n-1}\,I\,\left[1 + o(1)\right].
\end{equation}   
   
\noindent Finally, Equations \eqref{presquemaj} and \eqref{mindpsipsi} allow to   
estimate $J(\bar\psi_\ep)$ in the following way:   
   
$$J(\psi_\ep)=\frac{\Big(\int_M|\bar   
  D\bar\psi_\ep|^{\frac{2n}{n+1}}\mathrm{d}v_g\Big)^\frac{n+1}{n}}{\int_M\Re e<\bar   
  D\bar\psi_\ep,\bar\psi_\ep>\mathrm{d}v_g} \leq 
\frac{n}{2}\,\omega_{n-1}^\frac{1}{n}\,I^\frac{1}{n}\left[1+o(1)\right]\;.$$   
   
\noindent By (\ref{conform}), we have  
$$w_{n-1} I = \int_{\mR^n} f^n dx= \omega_n$$

\noindent Therefore, we proved that for the test spinor $\varphi$, we have   
\begin{equation}   
  \label{finale}   
  J(\bar\psi_\ep) \leq \lamin(\mS^n)\,\left[1+o(1)\right]\;.   
\end{equation}

\noindent Hence Theorem~\ref{main} is proven.   
 
\begin{rem}
There is a variant of this proof which needs less calculations.
As a first step, one proves that for any $\ep>0$ there is a test spinor
$\phi_\ep$ on $\mR^n$ with support in $B_0(1)$ such that 
$J^{\mR^n}_{\geucl}(\phi_\ep)\leq \lamin(\mS^n)+ \ep$ where $\ep>0$. The argument
for this coincides with the above proof, but the terms $\bf IV$ to $\bf XXI$
vanish, as $\mR^n$ is flat.

In a second step, one transplants this compactly supported spinor $\phi_\ep$
to the arbitrary compact spin manifold $(M,\Lambda^2 g)$, 
where $\Lambda>0$ is constant, and one obtains a spinor $\ol{\phi_\ep}$ on
$(M,\Lambda^2 g)$.
The terms $\bf IV$ to $\bf XXI$ reappear. 
However, from our Taylor expansion worked before, it is easy to see
that for $\Lambda \to \infty$ these terms dissapear.

One concludes that there for any $\ep>0$ there is a $\Lambda_\ep>0$ and
a spinor $\ol{\phi_\ep}$ on $(M,\Lambda_\ep^2 g)$ such that 
  $$J_{\Lambda_\ep^2g}(\ol{\phi_\ep})<\lamin(\mS^n) + 2 \ep.$$
Together with
  $$\lamin(M,g,\si)=\lamin(M,\Lambda_\ep^2g,\si)\leq 
    J_{\Lambda_\ep^2g}(\ol{\phi_\ep})$$
the theorem follows.

This proof is simpler. We chose the way presented above 
because of various reasons. One the other hand, 
as indicated in the introduction, 
in the case $n\geq 3$ it is not the result, but the method of proof which is 
interesting. The above formulae enter at several places 
in the literature, e.g.\  
\cite{ammann.humbert.morel:p03b}, \cite{ammann.humbert:03} and
\cite{raulot:06}.
Secondly, the simpler proof is close to Gro\ss{}e's 
proof \cite{grosse:06} and we refer to her article 
for the probably most elegant proof in dimension $n\geq 3$. Also in her proof 
some Taylor expansions from the present article are used.
\end{rem}

\section{The case $n=2$}  \label{dim2}
The 2-dimensional case is simpler since $g$ is locally conformally flat. On the other hand, 
some estimates of the last section are no longer valid in dimension~$2$, hence some
parts have to be modified. These modifications will be carefully carried out in this section.

Let
$(M,g)$ be a compact Riemannian surface equipped with a spin structure. 
If $\bar{g}$ is conformal to $g$ we
denote by $\mu_1(\bar{g})$  the smallest
positive eigenvalue of $\Delta_{\bar{g}}$. 
We prove the theorem.
\begin{thm} \label{main2}
There exists a family of  metrics $(g_{\ep})_\ep$ conformal to $g$ for which 
\begin{eqnarray*} \
\limsup_{\ep \to 0} \la^+_1(g_{\ep})^2  \Vol_{g_{\ep}}(M)
\leqslant     4  \pi
\end{eqnarray*} 
\begin{eqnarray*}
 \liminf_{\ep \to 0} \ \mu_1(g_{\ep})\Vol_{g_{\ep}}(M)  \geqslant   8 \pi.
\end{eqnarray*}
\end{thm}

Theorem \ref{main2} clearly implies Theorem \ref{main}. \\
Roughly, these metrics can be described as follows. At first we choose
a metric in the conformal class which is flat in a neighborhood of a point $p$.
We remove a small ball around it and glue in a large truncated sphere. 
This removal and gluing can be done in such a way that we stay within 
a conformal class. $\ep\to 0$. In the limit this truncated sphere is 
getting larger and larger compared to the original part of $M$.

Agricola, Ammann and Friedrich asked
the following question \cite{aaf:99}: \\

\noindent {\it Let $M$ be a two-dimensional torus equipped with a trivial spin structure,
can we find on $M$ a Riemannian metric $\tilde{g}$
for which $\la^+_1(\tilde{g})^2 < \mu_1(\tilde{g})$? }

To understand this question, 
recall that the two-dimensional torus carries 4 spin structures. 
Three of them (the non-trivial ones) are spin boundaries: for these spin structures 
it is easy to find flat examples with 
$\la^+_1(\tilde{g})^2 =\frac{1}{4} \mu_1(\tilde{g})$. For the trivial spin structure, 
one has $\la^+_1(\tilde{g})^2 = \mu_1(\tilde{g})$ for all flat metrics and 
$\la^+_1(\tilde{g})^2 > \mu_1(\tilde{g})$ for many $S^1$-equivariant one's.

Clearly, Theorem \ref{main2} answers this question 
but says  much more: firstly, the result is true on any
compact Riemannian surface equipped with a spin structure and not only when
$M$ is a two-dimensional torus. In addition, the metric $\tilde{g}$  can be
chosen in a given conformal class. Finally, this metric $\tilde{g}$  can be
chosen such that 
$(2-\de) \la^+_1(g)^2 < \mu_1(g)$ where $\de>0$ is arbitrarily
small. More precisely Theorem \ref{main} shows the corollary 

\begin{cor}[Proposition~\ref{prop.compar} of the Introduction]\label{compar}
On any compact Riemannian surface $(M,g)$, we have 
$$\inf \frac{\la^+_1(\bar{g})^2 }{\mu_1(\bar{g})} \leqslant  \frac{1}{2}$$
where the infimum is taken over all metrics $\bar{g}$ conformal to $g$.
\end{cor}

\subsection{$C^0$-metrics}
Let $f$ be a smooth positive function and set $\bar{g} = f^2 g$. Let also 
for $u \in C^{\infty}(M)$ 
$$I_{\bar{g}}(u) = \frac{\int_M | \nabla u |_{\bar{g}} dv_{\bar{g}} }{\int_M u^2
  dv_{\bar{g}} }.$$ 
It is well known that 
$\mu_1(\bar{g}) = \inf I_{\bar{g}} (u)$ where the infimum is taken over the
smooth non-zero functions $u$ for which $\int_M u dv_{\bar{g}} = 0$. 
Another way to express $\mu_1(\bar{g})$ is
\begin{equation}\label{varcharmu}
\mu_1(\bar{g}):=\inf_V \sup_{u\in V\setminus\{0\}}  I_{\bar{g}} (u)
\end{equation}
where the infimum runs over all $2$-dimensional subspaces $V$ of $C^\infty(M)$.
We now can write  all these expressions in the metric $g$. We then see that
for  $u \in C^{\infty}(M)$, we have 
$$I_{\bar{g}}(u) = \frac{\int_M |\nabla u|_g^2 dv_g}{\int u^2 f^2 dv_g}$$
and $\mu_1(\bar{g})$ is characterized in a way analogous to \eref{varcharmu}.
Now if $f$ is only continuous, we can define  $\bar{g} = f^2 g$. The
symmetric $2$-tensor $\bar{g}$ is not really a metric since $f$ is not smooth. We then
say that $g$ is a {\it $C^0$-metric}. We can define the first eigenvalue
$\mu_1(\bar{g})$ of $\Delta_{\bar{g}}$ using the definition above. 

Suppose that 
\begin{equation}\label{tiffass}
(1+\rho)^{-1}f\leq {\ti f} \leq (1+\rho) f.
\end{equation}
Then  
$$(1+\rho)^{-2} I_{{\ti f}^2g}(u)\leq I_{{f}^2g}(u) \leq (1+\rho)^2 I_{{\ti f}^2g}(u).$$
\noindent From the variational characterization \eref{varcharmu} it the follows that 
  $$(1+\rho)^{-2}\mu_1({\ti f}^2g)\leq \mu_1(f^2g)\leq (1+\rho)^2\mu_1({\ti f}^2g),$$
which is a special case of a result by Dodziuk \cite[Proposition 3.3]{dodziuk:82}.
In particular, we get
\begin{lem} \label{gene_lapla}
If  $(f_n)$ is  a
sequence of smooth positive functions that converges uniformily to $f$,
then  $\mu_1(f_n^2 g)$ tends to $\mu_1(f^2g)$.
\end{lem} 


\noindent In the same way, if $\bar{g} = f^2 g$ is a  metric conformal to $g$
where $f$ is positive and smooth, we define 
$$\cJ_{\bar{g}}(\psi) = \frac{\displaystyle\int_M |D_{\bar{g}}\psi|_{\bar{g}}^2
  dv_{\bar{g}}}{\displaystyle\int_M \langle  D_{\bar{g}} 
  \psi, \psi\rangle _{\bar{g}} dv_{\bar{g}}}.$$
The first eigenvalue of the Dirac operator $D_{\bar{g}}$ is then given by 
$\lambda_1^+(\bar{g})  = \inf \cJ_{\bar{g}}(\psi)$ where the infimum is taken over the
smooth spinor fields $\psi$ for which $  \int_M \langle  D_{\bar{g}} \psi, \psi\rangle 
dv_g > 0$. 
Now, as explained in paragraph~\ref{subsec.confch} 
we can identify spinors for the
metric $g$ and spinors for the metric $\bar{g}$ by a fiberwise isometry. 
Moreover, using this identification, we have
for all smooth spinor field:
$$D_{\bar{g}}( f^{-\frac{1}2} \phi )= f^{ -\frac{3}2}D_g \phi.$$
This implies that if we set $\phi = f^{\frac{1}{2}} \psi$, we have   
$$\cJ'_{\bar{g}}(\phi): = \frac{\displaystyle\int_M |D_g\phi|^2 f^{-1} dv_g}{\displaystyle\int_M \langle  D_g \phi, \phi\rangle  dv_g}= \cJ_{\bar{g}}(\psi)$$
and 
the first eigenvalue of the Dirac operator $D_{\bar{g}}$ is given by 
\begin{equation}\label{def.laone}
  \lambda_1^+(\bar{g})  = \inf \cJ'_{\bar{g}}(\phi)
\end{equation} where the infimum is taken over the
smooth spinor fields $\phi$ for which $ \displaystyle 
\int_M \langle  D_g \phi, \phi\rangle 
dv_g>0$.
Now, when  $\bar{g}=f^2 g$ is no longer smooth, but a $C^0$-metric,
we can use \eref{def.laone} to define $\la^+_1(\bar{g})$.

Under the assumption \eref{tiffass}, we get
  $$(1+\rho)^{-1}\cJ'_{\ti f^2g}(\phi) \leq \cJ'_{f^2g}(\phi) \leq (1+\rho)\cJ'_{\ti f^2g}(\phi),$$
and hence 
  $$(1+\rho)^{-1}\lambda_1^+(\ti f^2g)\leq \lambda_1^+(f^2g)\leq (1+\rho)\lambda_1^+(\ti f^2g).$$
We have proven a result similar to Lemma \ref{gene_lapla}:

\begin{lem} \label{gene_dir}
If  $(f_n)$ is  a
sequence of smooth positive functions that converges uniformly to $f$,
then  $\la^+_1(f_n^2 g)$ tends to $\la^+_1(\bar{g})$.
\end{lem}



\subsection{The metrics $(g_{\al,\ep})_{\al,\ep}$}
In this paragraph, we construct the metrics $(g_{\al,\ep})_{\al,\ep}$
conformal to $g$ which
will satisfy:

\begin{eqnarray} \label{main_dir} 
\limsup_{\ep \to 0} \la^+_1(g_{\al,\ep})^2  \Vol_{g_{\al,\ep}}(M)
\leqslant    4 \pi 
\end{eqnarray} 
and 
\begin{eqnarray}\label{main_lapla}
\liminf_{\al \to 0} 
 \liminf_{\ep \to 0} \mu_1(g_{\al,\ep})\Vol_{g_{\al,\ep}}(M)  
\geqslant   8 \pi.
\end{eqnarray}
Clearly this implies 
Theorem \ref{main}.
By Lemmas \ref{gene_lapla} and \ref{gene_dir}, it suffices to construct $C^0$-metrics $(g_{\al,\ep})_{\al,\ep}$. Recall that the volume of $M$ for 
a $C^0$-metric is defined by
$\Vol_{f^2 g} (M)= \int_M f^2 dv_g$. 
At first, without loss of generality, we can assume that $g$ is
flat near a point $p \in M$. Let $\al>0$ be a small number to be fixed
later such that $g$
is flat on $B_p( \al)$.  We set  
for all $x \in M$ and $\ep >0$,
\[ f_{\al,\ep}(x) = \left\{\begin{array}{ccc}
\frac{\ep^2}{\ep^2 + r^2} & \text{if} & r \leqslant   \al\\
\frac{\ep^2}{\ep^2 + \al^2} & \text{if} & r >  \al
\end{array} \right. \]
where $r = d_g(.,p)$. The function $f_{\al,\ep}$ is continuous and 
positive on $M$. We then define for all $\ep >0$, 
$g_{\al,\ep} = f_{\al,\ep}^2 g$. 
The symmetric 2-tensors $(g_{\al,\ep})_{\al,\ep}$ will be the desired
$C^0$-metrics.
For these metrics, we have 
$$\Vol_{g_{\al,\ep}}(M) = \int_M f_{\al,\ep}^2dv_g = \int_{B_p(\al)} f_{\al,\ep}^2 dv_g 
+ \int_{M \setminus B_p(\al)} f_{\al,\ep}^2 dv_g.$$
Since $g$ is flat on $B_p(\al)$, we have 
$$\int_{B_p(\al)} f_{\al,\ep}^2 dv_g = \int_0^{2 \pi} \int_0^{\al}
\frac{\ep^4 r}{(\ep^2 + r^2)^2} dr d\Theta.$$
Substituting $\rho=r^2/\ep^2$ we obtain
$$\int_{B_p(\al)} f_{\al,\ep}^2 dv_g  = \pi \ep^2  \int_0^{\frac{\al^2}{\ep^2}}
\frac{1}{(1 + \rho)^2} dr =  \pi \ep^2 + o(\ep^2).$$ 
Since  $f_{\al,\ep}^2  \leqslant  \frac{\ep^4}{\al^4}$ on $M \setminus B_p(\al)$,
 we have $\int_{M \setminus B_p(\al)}
f_{\al,\ep}^2 dv_g = o(\ep^2)$. 
We obtain 
\begin{eqnarray} \label{volume}
 \Vol_{g_{\ep}}(M)=\pi \ep^2 + o(\ep^2).  
\end{eqnarray}

\subsection{Proof of relation (\ref{main_dir})}
We define on $\mR^2$ as in subsection~\ref{subsec.explicit} the spinor field
$$\psi(x) = f(x) (1 - x) \cdot \psi_0$$
where $f(x)=\frac{2}{1+|x|^2}$, $|\psi_0|=1$. We have 
\begin{eqnarray} \label{test_spinor} 
 D\psi  = f \psi \; \hbox{ and } \;  |\psi|= f^{\frac{1}{2}}.
\end{eqnarray}
Now, we fix a small number $\al >0$ such that $g$ is flat on
$B_p(2 \al)$. Then, let $\de$ be a small number such that 
we take $0 \leqslant  \de \leqslant \alpha$. Assume that $\ep$
tends to $0$. Furthermore let 
$\eta$ be a smooth cut-off function defined on $M$ by 
\[ \eta (x) = \left| \begin{array}{ccc}
 1 & \mbox{if} & r \leq \de^2\hfill \\
 \frac{\log(r)}{\log(\de)}-1 & \mbox{if} & r \in [\de^2, \de] \\
0 & \mbox{if} & r \geq \de\hfill 
\end{array} \right. \]
The function $\eta$
is such that 
$0 \leqslant  \eta \leqslant  1$, $\eta(B_p(\de)) = \{1 \}$, $\eta(\mR^n \setminus B_p(2 \de))= \{0\}$ and
\begin{eqnarray} \label{gradto0}
\kappa_\delta:=  \int_{M} |\nabla \eta|^2 dv_g\to 0 \mbox{ for }\delta\to 0.
\end{eqnarray}
 Identifying $B_p(2\de)$ in $M$ with $B_0(2\de)$ in
  $\mR^2$, we can define a smooth spinor field on $M$ by $\psi_{\ep} =
  \eta(x) \psi\left( \frac{x}{\ep}\right)$. 
Using (\ref{test_spinor}), we have 
\begin{eqnarray} \label{dpsi}
D_g(\psi_{\ep}) = \nabla \eta \cdot \psi\left(\frac{x}{\ep}\right) + \frac{\eta}{\ep}
f\left(\frac{x}{\ep}\right) \psi\left(\frac{x}{\ep}\right).
\end{eqnarray}
Since $\langle  \nabla \eta \cdot \psi(\frac{x}{\ep}),  \psi(\frac{x}{\ep})\rangle  \in
i \mR$ and since $|D_g \psi_{\ep}|^2 \in \mR$, we have 
\begin{eqnarray} \label{intdpsi}
\int_M |D_g \psi_{\ep}|^2 f_{\al,\ep}^{-1} dv_g = I_1 + I_2 
\end{eqnarray} 
where 
$$I_1 = \int_M |\nabla \eta|^2 \left| \psi\left(\frac{x}{\ep}\right)\right|^2 dx \; \hbox{ and }
I_2 = \int_M \frac{\eta^2}{\ep^2} f^2 \left(\frac{x}{\ep}\right) \left|
\psi\left(\frac{x}{\ep}\right)\right|^2 f_{\al,\ep}^{-1} dx.$$
By (\ref{test_spinor}), $\left| \psi\left(\frac{x}{\ep}\right)\right|^2
\leq 2$ and hence  
\begin{eqnarray} \label{I1}
I_1 \leqslant  2 \int_M |\nabla \eta|^2 dv_g = 2\kappa_\de\to 0
\end{eqnarray}
for $\de \to 0$.
Now, by (\ref{test_spinor}),
$$I_2 \leqslant  \frac{2}{\ep^2} \int_{B_p(2\de)}  f^3 (\frac{x}{\ep})f_{\al,\ep}^{-1} dx.$$
Since $f_{\al,\ep} =   \frac{1}{2} f(\frac{x}{\ep})$ on the support of
$\eta$,  we have 
$$I_2 \leqslant  \frac{2}{\ep^2} \int_{B_p(2\de)} f^2 (\frac{x}{\ep}) dx.$$
Mimicking what we did to get (\ref{volume}), we obtain that 
$$I_2 \leqslant  8 \pi + o_\ep(1)$$
where $o_\ep(1)$ denotes a term tending to $0$ for $\ep\to 0$.
Together with (\ref{intdpsi}) and (\ref{I1}), we obtain  
\begin{eqnarray} \label{intdpsi2}
\int_M |D_g \psi_{\ep}|^2 f_{\al,\ep}^{-1} dv_g \leqslant 8 \pi + 2\kappa_\de + o_\ep(1).
\end{eqnarray}
In the same way, by (\ref{dpsi}), since $\displaystyle\int_M \langle D_g (\psi_\ep), \psi_\ep\rangle 
dv_g \in \mR$ and since $\langle  \nabla \eta \cdot \psi(\frac{x}{\ep}),  \psi(\frac{x}{\ep})\rangle  \in
i \mR$, we have 
$$\int_M \langle D_g (\psi_\ep), \psi_\ep\rangle 
dv_g  =  \int_M \frac{\eta^2}{\ep}   f\left(\frac{x}{\ep}\right) \left|
\psi\left(\frac{x}{\ep}\right)\right|^2 dv_g .$$
By (\ref{test_spinor}), this gives 
$$\int_M \langle D_g (\psi_\ep), \psi_\ep\rangle 
dv_g  = \int_M \frac{\eta^2}{\ep} f^2\left(\frac{x}{\ep}\right)dv_g.$$
With the computations made above, it follows that 
$$\int_M \langle D_g (\psi_\ep), \psi_\ep\rangle 
dv_g = 4 \pi \ep + o(\ep).$$ 
Together with (\ref{intdpsi2}) and (\ref{volume}), we obtain
\begin{eqnarray*}
\la^+_1(g_{\al,\psi})^2  \Vol_{g_{\al,\psi}}(M) &\leqslant&  
\left(\cJ'_{g_{\al,\ep}}(\psi_{\ep})\right)^2 \Vol_{g_{\al,\ep}}(M)  \leqslant 
  \left( \frac{8 \pi + 2\kappa_\de + o_\ep(1)}
   {4 \pi \ep+ o(\ep)} \right)^2 (\pi \ep^2 + o(\ep^2))\\
  &=& 4 \pi + 2\kappa_\de + \frac{1}{4\pi} \kappa_\delta^2 +o_\ep(1).
\end{eqnarray*}
Letting $\ep$ then $\de$ go to $0$, we get Relation (\ref{main_dir}).

\subsection{Proof of relation (\ref{main_lapla})} As pointed out by 
the referee the metrics $g_{\al,\ep}$ coincide with metrics 
contructed in \cite{takahashi:02}, and 
relation (\ref{main_lapla}) is proven in this article.

\vspace{1cm}             
Authors' addresses:             
\nopagebreak 
\vspace{5mm}\\ 
\parskip0ex 
\vtop{ 
\hsize=11cm\noindent 
\obeylines             
Bernd Ammann, Jean-Fran\,cois Grosjean, and Emmanuel Humbert,             
Institut \'Elie Cartan BP 239            
Universit\'e de Nancy 1             
54506 Vandoeuvre-l\`es -Nancy Cedex             
France                   
}      
\vtop{ 
\hsize=10cm\noindent 
\obeylines             
Bertrand Morel
D\'epartements QLIO et GIM
IUT Epinal - Hubert Curien
7, rue des Fusill\'es de la R\'esistance
BP 392
88010 Epinal Cedex
France
}
\vspace{0.5cm}             
\nopagebreak 
            
E-Mail:             
{\tt ammann@iecn.u-nancy.fr}, {\tt grosjean@iecn.u-nancy.fr}, 
{\tt humbert@iecn.u-nancy.fr}, \\
and {\tt bertrand.morel@univ-nancy2.fr} 
                  

\begin{thebibliography}{AHM60}           
         
\bibitem[AAF99]{aaf:99}
I. Agricola, B. Ammann and T. Friedrich,
{A comparison of the eigenvalues of the Dirac and Laplace operators on
  a two-dimensional torus},
\newblock{\em Manuscripta Math.} \textbf{100} (1999), No 2,  231--258.

\bibitem[Am03a]{ammann.habil}   
B. Ammann, \emph{A variational problem in conformal spin Geometry}, 
Habilitationsschrift, Universit\"at Hamburg, May 2003,  
downloadable on http://www.berndammann.de/publications. 

\bibitem[Am03b]{ammann:03}   
B. Ammann, {A spin-conformal lower bound of the first positive {D}irac   
  eigenvalue}, \newblock{\em {D}iff. {G}eom. {A}ppl.} \textbf{18} (2003), 21--32.   
\bibitem[Am03c]{ammann:p03}
B. Ammann, {The smallest {D}irac eigenvalue in a spin-conformal class
and cmc-immersions}, {{P}reprint}, {2003}.

\bibitem[AH03]{ammann.humbert:03}         
B. Ammann, E. Humbert,   
\newblock{Positive mass theorem for the Yamabe problem on spin
  manifolds,} \newblock{\em Geom. Funct. Anal.} \textbf{15} (2005),
567--576.       

\bibitem[AHM03]{ammann.humbert.morel:p03b}         
B. Ammann, E. Humbert, B. Morel,  
\newblock{Mass endomorphism and spinorial Yamabe type problems  
on conformally flat manifolds,} 
\newblock{\em Comm. Anal. Geom.} \textbf{14} (2006),  163--182. 
     
\bibitem[Aub76]{aubin:76}           
T. Aubin,           
\newblock {{\'E}quations diff{\'e}rentielles non lin{\'e}aires et probl{\`e}me           
  de Yamabe concernant la courbure scalaire,}           
\newblock {\em J. Math. Pur. Appl., IX. Ser.} \textbf{55} (1976), 269--296. 
   
\bibitem[B\"ar91]{baer:thesis}           
C. B\"ar, {{D}as {S}pektrum des {D}irac-{O}perators}, PhD thesis,
{\em Bonner Math. Schr.},
{\bf 217}, Bonn 1991.

\bibitem[B\"ar92]{baer:92b}           
C. B\"ar,           
\newblock {Lower eigenvalue estimates for Dirac operators,}           
\newblock {\em Math. Ann.} {\textbf 293} (1992), 39--46.           
    
\bibitem[BG92]{bourguignon.gauduchon:92} 
J.-P. Bourguignon and P. Gauduchon, 
\newblock {Spineurs, op{\'e}rateurs de {D}irac et variations de m{\'e}triques,} 
\newblock {\em Commun. Math. Phys.} \textbf{144} (1992), 581--599. 
  
\bibitem[CES03]{colbois.elsoufi:03} 
B. Colbois and A. El Soufi,
\newblock{Extremal eigenvalues of the Laplacian in a conformal class of
  metrics: the 'conformal spectrum'}
\newblock{\em Ann. Global Anal. Geom.} \textbf{24} (2003),  337--349.  

\bibitem[Dod82]{dodziuk:82}
J.~Dodziuk, {Eigenvalues of the {L}aplacian on forms},
{\em Proc. Amer. Math. Soc.} \textbf{85} (1982), 437--443.

\bibitem[Fri80]{friedrich:80}
T.~Friedrich, {Der erste {E}igenwert des {D}irac-{O}perators
einer kompakten {R}iemannschen {M}annigfaltigkeit nicht-negativer {K}r{\"u}mmung}, {\em Math. Nach.},
{\bf 97} (1980),
117--146.

\bibitem[Gro06]{grosse:06}
N. Gro\ss{}e, {On a conformal invariant of the Dirac operator on noncompact manifolds,}
{\em Ann. Glob. Anal. Geom.} \textbf{30} (2006), 407--416.

\bibitem[Heb97]{hebey:97}           
E. Hebey,           
\newblock  {\em Introduction \`a l'analyse non-lin\'eaire sur les vari\'et\'es,}   
\newblock {Diderot \'Editeur, Arts et sciences.}, 1997.        
   
\bibitem[Hij86]{hijazi:86}           
O.~Hijazi,           
  \newblock  A conformal lower bound for the smallest eigenvalue of the {D}irac operator and {K}illing Spinors,            
\newblock {\em Commun. Math. Phys.} \textbf{104} (1986), 151--162. 
           
\bibitem[Hij91]{hijazi:91}             
O.~Hijazi,             
\newblock  Premi\`ere valeur propre de l'op\'erateur de {D}irac et nombre de {Y}amabe,             
\newblock {\em C. R. Acad. Sci. Paris}, \newblock{\em S\'erie I } \textbf{313}, (1991), 865--868.            
 
\bibitem[Hij01]{hijazi:01}
O.~Hijazi,              
\newblock Spectral properties of the {D}irac operator and geometrical structures,
\newblock{\em Ocampo, Hernan (ed.) et al., Geometric methods for quantum field theory}, 
Proceedings of the summer school, Villa de Leyva, Colombia, July 12-30, 1999. 
Singapore: World Scientific. 116-169,
2001.

\bibitem[Hit74]{hitchin:74}
N.~Hitchin, 
\newblock Harmonic spinors, 
\newblock{\em Adv. Math.} \textbf{14} (1974), 1--55.

\bibitem[LM89]{lawson.michelsohn:89}
H.-B. Lawson and M.-L. Michelsohn,
{\em Spin Geometry},
{Princeton University Press},
{Princeton},
{1989}.

\bibitem[LP87]{lee.parker:87}           
J.~M. Lee and T.~H. Parker,           
\newblock {The Yamabe problem,}           
\newblock {\em Bull. Am. Math. Soc., New Ser.} \textbf{17} (1987), 37--91.           
 
\bibitem[Lot86]{lott:86}           
J. Lott,           
\newblock {Eigenvalue bounds for the Dirac operator,}           
\newblock {\em Pacific J. Math.} \textbf{125} (1986), 117--126.             

\bibitem[Pfa02]{pfaeffle:thesis}
F. Pf{\"a}ffle,
{\em {E}igenwertkonvergenz f{\"u}r {D}irac-{O}peratoren},
PhD thesis,
University of {H}amburg, {G}ermany, 2002,
Shaker Verlag Aachen 2003, ISBN 3-8322-1294-9.

\bibitem[Rau06]{raulot:06}
S. Raulot, 
\newblock {Aspect conforme de l'op\'erateur de Dirac sur une vari\'et\'e 
\`a bord}, PhD thesis, Universit\'e Nancy I, 2006.  
            
\bibitem[Sch84]{schoen:84}           
R. Schoen,           
\newblock {Conformal deformation of a Riemannian metric to constant scalar           
  curvature},           
\newblock {\em J. Diff. Geom.} \textbf{20} (1984), 479--495.           
           
\bibitem[Tru68]{trudinger:68}           
N.S. Trudinger,           
\newblock {Remarks concerning the conformal deformation of Riemannian           
  structures on compact manifolds},  
\newblock {\em Ann. Sc. Norm. Super. Pisa, Sci. Fis. Mat., III. Ser.}           
  \textbf{22} (1968), 265--274.           
           
\bibitem[Tak02]{takahashi:02}           
J. Takahashi,           
\newblock {Collapsing of connected sums and the eigenvalues of the
  Laplacian}, \newblock {\em J. Geom. Phys.} \textbf{40} (2002),
201--208.           
           
\bibitem[Yam60]{yamabe:60}           
H. Yamabe,           
\newblock {On a deformation of Riemannian structures on compact manifolds,}           
\newblock {\em Osaka Math. J.} \textbf{12} (1960), 21--37.                     
\end{thebibliography}
\end{document}